\documentclass[10pt,twoside]{article}
\usepackage{amsmath,amsthm}
\usepackage{picinpar}
\usepackage{longtable}
\usepackage[dvips]{graphicx}
\usepackage[mathscr]{eucal}
\usepackage{calc}
\usepackage{ifthen}
\usepackage{enumerate}
\usepackage{amssymb,latexsym}
\usepackage{amsthm}
\usepackage[dvips]{graphics}
\usepackage{color}
\usepackage{tabularx}
\usepackage{makeidx}

\theoremstyle{definition}

 \newtheorem{teor}{Theorem}

\newtheorem{corol}{Corollary}

\newtheorem{Nota}{Remark}

\begin{document}
\thispagestyle{plain}
\par\bigskip
\begin{centering}

\textbf{ON SUMS OF FIGURATE NUMBERS BY USING TECHNIQUES OF POSET REPRESENTATION THEORY}

\end{centering}\par\bigskip
 \begin{centering}
 \footnotesize{AGUSTIN MORENO CA\~NADAS}\par\bigskip
 \end{centering}

$\centerline{\textit{\small{Departamento de Matem\'aticas, Universidad Nacional de Colombia}}}$
\centerline{\textit{\small{Bogot\'a-Colombia}}} $\centerline{\textit{\small{amorenoca@unal.edu.co}}}$

\bigskip

\bigskip \small{We use representations and differentiation algorithms of posets, in order to obtain
results concerning unsolved problems on figurate numbers. In particular, we present criteria for natural numbers
which are the sum of three octahedral numbers, three polygonal numbers of positive rank or four cubes with two of
them equal. Some identities of the Rogers-Ramanujan type involving this class of numbers are also obtained.}


\bigskip

\textbf{1. Introduction}\par\bigskip

The representation theory of posets was developed in the 70's by Nazarova and Roiter [16]. The main aim of their
study was to determine indecomposable matrix representations of a given poset $\mathscr{P}$ over a fixed field
$k$.\par\bigskip

An important role in poset representation theory is played by some differentiation algorithms, reducing the study
of matrix representations of a given dimension type to representations of smaller dimension type. One of the
algorithms is the algorithm of differentiation with respect to a pair of points D-I due to Zavadskij which
associates to any poset $\mathscr{P}$ with a suitable pair $(a,b)$ of elements $a,b\in\mathscr{P}$ a poset
$\mathscr{P}'_{(a,b)}$. The algorithms are successfully applied in determining the representation type (finite or
infinite) of posets and in the classification of indecomposable poset representations.\par\bigskip Soon after the
discovery of matrix representations of a poset $\mathscr{P}$, Gabriel introduced the concept of a filtered
$k$-linear representation of $\mathscr{P}$ in connection with the investigation of oriented graphs having finitely
many isomorphism classes of indecomposable linear representations. In this case a $k$-linear representation is a
system $U=(U_{0}\hspace{0.1cm};\hspace{0.1cm}U_{x}\mid x\in\mathscr{P})$, where $U_{0}$ is a finite dimensional
$k$-vector space and for each $x\in\mathscr{P}$, $U_{x}$ is a subspace $U_{x}\subset U_{0}$ such that
$U_{x}\subset U_{y}$ if $x\leq y$ in $\mathscr{P}$. These two concepts of poset representations can be connected
by defining the associated matrix of a  $k$-linear representation.\par\bigskip

The purpose of this paper is to link the poset representation theory to the additive number theory, by adapting
the tools mentioned above in order to solve different types of diophantine equations which involve figurate
numbers, where a \textit{figurate number} is a number that can be represented by a regular geometrical arrangement
of equally spaced points. If the arrangement forms a regular polygon the number is called a \textit{polygonal
number} [2,5].\par\bigskip

The $n$-th polygonal number of order or rank $k$, $p^{n}_{k}$ is given by the formula (often, $0$ is included as a
polygonal number)
\par\bigskip

\begin{centering}
$p^{n}_{k}=\frac{1}{2}[(n-2)k^{2}-(n-4)k]$.\par\bigskip

\end{centering}

Perhaps the most remarkable result concerning polygonal numbers was stated by Fermat as early as 1638, who made
the statement that every number is expressible as the sum of at most $t$ $t$-gonal numbers, $t\geq3$ [5,9]. A
Fermat's proof of this fact has not been found yet. Meanwhile for triangular numbers, Gauss proved (1796), that
every number is expressible as the sum of three triangular numbers (i.e., every number $m$ is expressible as the
sum of three numbers of the form $\frac{k(k+1)}{2}$). Gauss's statement is equivalent to the statement that every
number of the form $8m+3$ can be expressible as the sum of three odd squares. About this particular fact, Legendre
(1798) published first that a number of the form $4^{k}(8m+7)$, $k,m\geq0$ cannot be expressible as the sum of
three squares and the proof of this result was given by Gauss (1801) in his Disquisitiones Arithmeticae
[5].\par\bigskip

Euler (1772) stated that at least $a+2n-2$ terms are necessary to express every number as a sum of figurate
numbers
\par\bigskip

\begin{centering}
$1$, $n+a$, $\frac{(n+1)(n+2a)}{1\cdot2}$, $\frac{(n+1)(n+2)(n+3a)}{1\cdot2\cdot3}$, $\dots$\par\bigskip
\end{centering}

Furthermore, Euler proved that if $n=\frac{k(k+1)}{2}$ is a triangular number then $9n+1$, $25n+3$, $49n+6$ and
$81n+10$ are also triangular numbers [7]. In fact, for every odd number $m=2j+1$ the number
$m^{2}n+\frac{m^{2}-1}{8}$ is a triangular number and each sum of the form
$n=\frac{j(j+1)}{2}+\frac{k(k+1)}{2}=p^{3}_{j}+p^{3}_{k}$ gives raise to an infinite sequence of this type of sums
for an odd number $m$, in such a way that\par\bigskip

\begin{centering}
$m^{2}n+\frac{m^{2}-1}{4}=(m^{2}p^{3}_{j}+\frac{m^{2}-1}{8})+(m^{2}p^{3}_{k}+\frac{m^{2}-1}{8})$.\par\bigskip

\end{centering}

Moreover, $n$ is the sum  of two triangular numbers exactly when\par\bigskip

\begin{centering}
$2(4n+1)=(2j+1)^{2}+(2k+1)^{2}$.\par\bigskip

\end{centering}

That every number is a sum of four squares was proved by Lagrange (1772) with the help of Euler's work. Finally
Cauchy (1813) gave the first proof of the Fermat's assertion on polygonal numbers [5]. \par\bigskip

On the number of representations of a number $n$ as the sum of polygonal numbers, Dirichlet gave a formula for the
number of ways in which $n$ can be expressed as the sum of three triangular numbers, and Jacobi gave the formulas
for two or four squares [5,9,10].\par\bigskip

Recently, Ewell (1992) stated that the number of representations of $n$ as the sum of two triangular numbers is
$d_{1}(4n+1)-d_{3}(4n+1)$ (where $d_{i}(n)$ is the number of divisors of $n$ which are congruent to $i$) [6],
Conway and Schneeberger (1993) proved (the fifteen theorem) that if a positive integer-matrix quadratic form
represents each of 1, 2, 3, 5, 6, 7, 10, 14, 15, then it represents all positive integers [5], Bhargava and Hanke
(2005) proved the 290-theorem (see [5]) and, Farkas (2006) proved that every positive integer can be written as
the sum of two squares plus one triangular number and that every positive integer can be written as the sum of two
triangular numbers plus one square. Furthermore Farkas, gave formulas for the number of ways in which $n$ can be
expressed as the sum of three triangular numbers or three squares (by using theta functions) [8].\par\bigskip
Although the Fermat's theorem has been proved, there are still open problems concerning this theorem and others
figurate numbers. In particular, we will investigate the following problems proposed by Guy in
[9,10,11]:\par\bigskip
\begin{enumerate}[(1)]
\item What theorems are there, stating that all numbers of a suitable shape are expressible as the sum of three
(say) squares of numbers of a given shape? For instance, can all sufficiently large numbers be expressed as the
sum of three pentagonal (hexagonal, heptagonal) numbers of nonnegative rank? Equivalently, is every sufficiently
large number of shape $24n+3$ $(8n+3,\hspace{0.1cm}40n+27)$ expressible as the sum of three squares of numbers of
shape $6r-1$ $(4r-1,\hspace{0.1cm}(10r\pm3))$?\par\bigskip

\item There are theorems giving the number of representations of a number $n$, as the sum of triangular or
square numbers. Can we find corresponding results for any of the other polygonal numbers?\par\bigskip

\item The Pollock octahedral numbers conjecture claims that every number is the sum of at most seven octahedral
numbers, where the $n$-th \textit{octahedral number} $\mathscr{O}_{n}$ is given by the formula
$\frac{n(2n^{2}+1)}{3}$. The Corresponding conjecture for tetrahedral numbers claims that every number is the sum
of at most five tetrahedral numbers, where the $n$-th \textit{tetrahedral number} $\rho_{n}$ is given by the
formula $\frac{n(n+1)(n+2)}{6}$. In this case Chou and Deng believe that all numbers greatest than 343867 are
expressible as the sum of four tetrahedral numbers [11].\par\bigskip
\item Is every number of the form $9n\pm4$ the sum of four cubes? Deshouillers, Hennecart, Landreau and Purnaba
believe that 7373170279850 is the largest integer which cannot be expressed as the sum of four nonnegative
integral cubes. Actually more demanding is to ask if every number is the sum of four cubes with two of them equal
[4,11,12,13].

\end{enumerate}
\par\bigskip
\textbf{2. Preliminaries}

\par\bigskip

\textbf{2.1. Posets}\par\bigskip

An \textit{ordered set} (or \textit{partially ordered set} or \textit{poset}) is an ordered pair of the form
$(\mathscr{P},\leq)$ of a set $\mathscr{P}$ and a binary relation $\leq$ contained in
$\mathscr{P}\times\mathscr{P}$, called the \textit{order} (or the \textit{partial order}) on $\mathscr{P}$, such
that $\leq$ is reflexive, antisymmetric and transitive [3,15]. The elements of $\mathscr{P}$ are called the
\textit{points} of the ordered set. In fact we shall assume in this work that $\mathscr{P}\neq\varnothing$ and
there is a bijective map (not necessary order-preserving) from a subset not empty of the set of all positive
integers $\mathbb{N}$, to the set of points of $\mathscr{P}$.\par\bigskip

An ordered set $C$ is called a \textit{chain} (or a \textit{totally ordered set} or a \textit{linearly ordered
set}) if and only if  for all $p,q\in C$ we have $p\leq q$ or $q\leq p$ (i.e., $p$ and $q$ are
comparable).\par\bigskip

Given an arbitrary point $a\in\mathscr{P}$, we define $a^{\triangledown}=\{x\in\mathscr{P}\mid a\leq x\}$,
$a_{\vartriangle}=\{x\in\mathscr{P}\mid x\leq a\}$. In general for $A\subset\mathscr{P}$,
$A^{\triangledown}=\underset{a\in A}\cup a^{\triangledown}$, $A_{\vartriangle}=\underset{a\in A}\cup
a_{\vartriangle}$.\par\bigskip

Let $\mathscr{P}$ be an ordered set and let $x,y\in\mathscr{P}$ we say $x$ is covered by $y$ if $x<y$ and $x\leq
z<y$ implies $z=x$.\par\bigskip

Let $\mathscr{P}$ be a finite ordered set. We can represent $\mathscr{P}$ by a configuration of circles
(representing the elements of $\mathscr{P}$) and interconnecting lines (indicating the covering relation). The
construction goes as follows.

\begin{enumerate}[(1)]
\item To each point $x\in\mathscr{P}$, associate a point $p(x)$ of the Euclidean plane $\mathbb{R}^{2}$, depicted
by a small circle with center at $p(x)$.

\item For each covering pair $x<y$ in $\mathscr{P}$, take a line segment $l(x,y)$ joining the circle  at $p(x)$ to
the circle at $p(y)$.

\item Carry out (1) and (2) in such a way that

\begin{enumerate}[(a)]
\item if $x<y$, then $p(x)$ is lower than $p(y)$,

\item the circle at $p(z)$ does not intersect the line segment $l(x,y)$ if $z\neq x$ and $z\neq y$.

\end{enumerate}

\end{enumerate}
A configuration satisfying (1)-(3) is called a \textit{Hasse diagram} or \textit{diagram} of $\mathscr{P}$. In the
other direction, a diagram may be used to define a finite ordered set; an example is given below.\par\bigskip

\begin{centering}
\begin{picture}(95,95)

 \setlength{\unitlength}{1pt}
  \setlength{\unitlength}{1pt}

\put(54,30){\line(0,1){20}} \put(54,27){\circle{5}} \put(54,56){\line(0,1){20}} \put(54,53){\circle{5}}
\put(54,79){\circle{5}} \put(52,26){\line(-1,-1){20}} \put(30,4){\circle{5}} \put(26,-12){$f$} \put(63,25){$c$}
\put(63,50){$d$} \put(63,78){$e$} \put(85,2){$b$} \put(111,-22){$a$} \put(78,4){\circle{5}}
\put(80,1){\line(1,-1){20}} \put(56,26){\line(1,-1){20}} \put(102,-21){\circle{5}}

\put(42,-35){Fig. 1}

\end{picture}

\par\bigskip
\end{centering}

\par\bigskip
\par\bigskip
\par\bigskip

 We have only defined diagrams for finite ordered sets. It is not possible to represent the whole of an infinite
ordered set by a diagram, but if its structure is sufficiently regular it can often be suggested
diagrammatically.\par\bigskip

Let $\mathscr{P}$ be a poset and $S\subset\mathscr{P}$. Then  \textit{$a\in S$ is a maximal element of $S$} if
$a\leq x$ and $x\in S$ imply $a=x$. We denote the set of maximal elements of $S$ by $\mathrm{Max}\hspace{0.1cm}S$.
If $S$ (with the order inherited from $\mathscr{P}$) has a top element, $\top_{S}$ (i.e., $s\leq\top_{S}$ for all
$s\in S$), then $\mathrm{Max}\hspace{0.1cm}S=\{\top_{S}\}$; in this case $\top_{S}$ is called the
\textit{greatest} (or \textit{maximum}) element of $S$, and we write $\top_{S}=\mathrm{max}\hspace{0.1cm}S$. A
\textit {minimal} element of $S\subset\mathscr{P}$ and $\mathrm{min}\hspace{0.1cm}S$, the \textit{least} (or
\textit{minimum}) element of $S$ (when these exist) are defined dually, that is reversing the order.\par\bigskip

Suppose that $\mathscr{P}_{1}$ and $\mathscr{P}_{2}$ are (disjoint) ordered sets. The \textit{disjoint union}
$\mathscr{P}_{1}+\mathscr{P}_{2}$ of $\mathscr{P}_{1}$ and $\mathscr{P}_{2}$ is the ordered set formed by defining
$x\leq y$ in $\mathscr{P}_{1}+\mathscr{P}_{2}$ if and only if either $x,y\in\mathscr{P}_{1}$ and $x\leq y$ in
$\mathscr{P}_{1}$ or $x,y\in\mathscr{P}_{2}$ and $x\leq y$ in $\mathscr{P}_{2}$. A diagram for
$\mathscr{P}_{1}+\mathscr{P}_{2}$ is formed by placing side by side diagrams for $\mathscr{P}_{1}$ and
$\mathscr{P}_{2}$.\par\bigskip

\textbf{2.2. Partitions}\par\bigskip

A \textit{partition} of a positive integer $n$ is a finite nonincreasing sequence of positive integers
$\lambda_{1},\lambda_{2},\dots,\lambda_{r}$ such that $\sum^{r}_{i=1}\lambda_{i}=n$. The $\lambda_{i}$ are called
the \textit{parts} of the partition [1]. A \textit{composition} is a partition in which the order of the summands
is considered.\par\bigskip

Often the partition $\lambda_{1},\lambda_{2},\dots,\lambda_{r}$ will be denoted by $\lambda$ and we sometimes
write $\lambda=(1^{f_{1}}2^{f_{2}}3^{f_{3}}\dots)$ where exactly $f_{i}$ of the $\lambda_{j}$ are equal to $i$.
Note that $\sum_{i\geq1}f_{i}i=n$.\par\bigskip

The \textit{partition function} $p(n)$ is the number of partitions of $n$. clearly $p(n)=0$ when $n$ is negative
and $p(0)=1$, where the empty sequence forms the only partition of zero. Let $H$ be a set of positive integers, we
denote $p(\underline{H},n)$ the number of partitions of $n$ that have all their parts in $H$, where
$\underline{H}$ is the set of all partitions whose parts lie in $H$. We let $\underline{H}(\leq d)$ denote the set
of all partitions in which no part appears more than $d$ times and each part is in $H$.\par\bigskip

Thus if $\mathbb{N}$ is the set of all positive integers then
$p(\underline{\mathbb{N}}(\leq1),n)=p(\mathscr{D},n)$ where $\mathscr{D}$ is the set of all partitions with
distinct parts. \par\bigskip Euler stated the following fact. If $H_{0}$ is the set of all odd positive integers
then $p(\underline{H_{0}},n)=p(\mathscr{D},n)$. In order to obtain the proof of this theorem is necessary to
consider the \textit{generating function} $f(q)$ for a sequence $a_{0}, a_{1}, a_{2}, a_{3},\dots$ defined as the
power series $f(q)=\sum_{n\geq0}a_{n}q^{n}$. Concerning this result Andrews asks for subsets of positive integers
$S,T$ such that $p(\underline{S},n)=p(\underline{T},n-1)$ [1].\par\bigskip

Given a partition $\lambda=(1^{f_{1}}2^{f_{2}}3^{f_{3}}\dots i^{f_{i}}\dots)$ with parts in a set $H$. In [14] it
was defined the \textit{derivative of substitution} of $\lambda$, $\lambda_{i}(u)=\frac{\partial\lambda}{\partial
i}(u)$ in such a way that if $u\in H$ then $\lambda_{i}(u)$ is a new partition with parts in $H$ such that
$\lambda_{i}(u)=(1^{f_{1}}2^{f_{2}}3^{f_{3}}\dots (i-1)^{f_{i-1}}ui^{f_{i}-1}(i+1)^{f_{i+1}}\dots)$. Furthermore,
we can make such a substitution with several parts at the same time. For instance if $u_{1},u_{2},\dots,u_{k}\in
H$ then \par\bigskip

\begin{centering}
 $\frac{\partial^{k}\lambda}{\partial
i_{1}\partial i_{2}\dots
\partial i_{k}}(u_{1},u_{2}\dots,u_{k})=(1^{f_{1}}\dots
u_{t}i_{t}^{f_{i_{t}}-1}\dots
 m^{f_{m}})$,\quad $1\leq t\leq k$. \par\bigskip
\end{centering}

and\par\bigskip
\begin{centering}
$\frac{\partial^{n}\lambda}{\partial i^{n}}(u)=(1^{f_{1}}\dots (i-1)^{f_{(i-1)}}ui^{f_{i}-n}(i+1)^{f_{(i+1)}}\dots
m^{f_{m}})$, $n\leq f_{i}$.\par\bigskip \end{centering}

In the sequel we will enunciate the famous Rogers- Ramanujan identities.

\begin{teor}
The partitions of an integer $n$ in which the difference between any two parts is at least 2 are equinumerous with
the partitions of $n$ into parts $\equiv1$ or $4\hspace{0.1cm}\mathrm{mod}\hspace{0.1cm}5$.

\end{teor}

\begin{teor}
The partitions of an integer $n$ in which each part exceeds 1 and the difference between any two parts is at least
2 are equinumerous with the partitions of $n$ into parts $\equiv2$ or
$3\hspace{0.1cm}\mathrm{mod}\hspace{0.1cm}5$.

\end{teor}

One of the main goals of the present paper, is to give some identities of these types for some classes of figurate
numbers.

\par\bigskip
\textbf{3. Partitions and Representations of posets over $\mathbb{N}$}\par\bigskip

Let $(\mathbb{N},\leq)$ be the set of all positive integers endowed with its natural order and
$(\mathscr{P},\leq')$ a poset. A \textit{representation} of $\mathscr{P}$ over $\mathbb{N}$ [14] is a system of
the form
\begin{equation}\label{representation}
\begin{split}
\Lambda&=(\Lambda_{0}\hspace{0.1cm};\hspace{0.1cm}(n_{x},\lambda_{x})\mid x\in\mathscr{P}),\\
\end{split}
\end{equation}

where $\Lambda_{0}\subset\mathbb{N}$, $\Lambda_{0}\neq\varnothing$, $n_{x}\in\mathbb{N}$,
$\lambda_{x}\in\underline{\Lambda_{0}}$ and $|\lambda_{x}|$ is the size of the partition $\lambda_{x}$. Further
\begin{equation}\label{conditions 1}
\begin{split}
x\leq'y&\Rightarrow n_{x}\leq
n_{y},\quad|\lambda_{x}|\leq|\lambda_{y}|,\quad\text{and}\quad\text{max}\{\lambda_{x}\}\leq\hspace{0.1cm}\text{max}\{\lambda_{y}\}.\\
\end{split}
\end{equation}

Let $(\mathscr{P},\leq')$ be a poset and $(a,b)$ be a pair of its incomparable points. If
$\mathscr{P}=a+b_{\vartriangle}+C$, where $C= \{c_{1}<' c_{2}\dots<' c_{n}\}\neq\varnothing$ is a chain and $a<'
c_{1}$, $b<'c_{1}$ then the pair $(a,b)$ is called \textit{$\mathrm{L}$-suitable} or \textit{suitable for
differentiation $\mathrm{L}$} [14], and the \textit{derivative poset} of the poset $\mathscr{P}$ with respect to
the pair $(a,b)$ is a poset\par\bigskip
\begin{centering}
$\mathscr{P}'_{(a,b)}=a^{-}+\mathscr{P}/(a^{\triangledown}) + C^{-}+C^{+}$,\par\bigskip
\end{centering}
where  $C^{-}=\{c^{-}_{1}<\dots<c^{-}_{n}\},\hspace{0.2cm}C^{+}=\{c^{+}_{1}<\dots<c^{+}_{n}\}$ are two new chains
(replacing the chain $C$) with the relations $a^{-}<c^{-}_{1}=a$,\hspace{0.2cm} $c^{-}_{i}<
c^{+}_{i}$,\hspace{0.2cm}$1\leq i\leq n$. It is assumed that each of the points $a^{-}$, $c^{-}_{i}$ ($c^{+}_{i}$)
inherits all order relations of the point $a$ ($c_{i}$), with the points of the poset $\mathscr{P}\backslash C$.
Fig. 2 provides a Hasse diagram of this differentiation.
\par\smallskip

 \setlength{\unitlength}{1pt}
  \setlength{\unitlength}{1pt}
\begin{centering}
\begin{picture}(95,95)

 \multiput(45,28)(18,18){1}{\circle{5}}
 \multiput(63,46)(18,18){1}{\circle{5}}
 \multiput(81,64)(18,18){1}{\circle{5}}

\multiput(101,82)(18,18){1}{\circle{5}}
  \multiput(47,29)(9,9){1}{\line(1,1){15}}
 \multiput(65,47)(9,9){1}{\line(1,1){15}}
\multiput(83,65)(9,9){1}{\line(1,1){16}}

          \put(33,24){$a$}
           \put(51,50){$c_{1}$}
            \put(91,88){$c_{n}$}
            \put(87,24){$b$}
          \put(170,-5){Fig. 2}
           \put(82,28){\circle{5}}
           \put(99,11){\line(-1,1){15}}
           \put(80,30){\line(-1,1){15}}
            \put(108,7){\circle{18}}
           \put(103,3){B}
           \put(170,45){$\xrightarrow[(a,b)]{\mathrm{L}}$}

 \multiput(298,34)(18,18){1}{\circle{5}}
 \put(306,32){$c^{-}_{1}$}
 \multiput(316,52)(18,18){1}{\circle{5}}
\put(341,68){$c^{-}_{n}$} \multiput(336,70)(18,18){1}{\circle{5}}

  \multiput(300,32)(9,9){1}{\line(1,-1){15}}
 \multiput(300,35)(9,9){1}{\line(1,1){15}}
\multiput(319,54)(9,9){1}{\line(1,1){15}}
           \put(317,16){\circle{5}}
           \put(322,16){$a^{-}$}
\multiput(277,52)(9,9){1}{\line(-1,-1){15}} \multiput(258,33)(9,9){1}{\line(-1,-1){15}}

\multiput(260,35)(18,18){1}{\circle{5}}
 \put(265,27){$b$}
 \multiput(278,54)(18,18){1}{\circle{5}}
 \put(264,54){$c^{+}_{1}$}
 \multiput(296,72)(18,18){1}{\circle{5}}

\multiput(316,90)(18,18){1}{\circle{5}} \put(300,90){$c^{+}_{n}$}

 \multiput(280,55)(9,9){1}{\line(1,1){15}}
\multiput(298,73)(9,9){1}{\line(1,1){16}} \multiput(296,36)(9,9){1}{\line(-1,1){16}}
 \multiput(334,71)(9,9){1}{\line(-1,1){17}}
\multiput(314,54)(9,9){1}{\line(-1,1){16}}
            \put(234,13){\circle{18}}
           \put(230,11){B}

\end{picture}
\end{centering}
\par\bigskip
\par\bigskip

In [14], it is assumed that if $\Lambda$ is a representation of a poset $\mathscr{P}$ with a pair $(a,b)$
$\mathrm{L}$-suitable, $t$ fixed, and $k_{x}\in\mathbb{N}$ for all $x\in\mathscr{P}$, then

\begin{equation}
\begin{split}
\lambda_{x}&=((k_{x})^{t}),\hspace{0.2cm}\text{for all}\hspace{0.1cm}x\in b^{\triangledown}+B,\\
\frac{\partial^{|\lambda_{a}|}\lambda_{x}}{\partial k^{|\lambda_{a}|}_{x}}(\lambda_{a})&=((\text{I}_{x})(\text{J}_{x})),\hspace{0.2cm}\text{if}\hspace{0.1cm}x\in C,\\
\text{I}_{x}&=\lambda_{a},\\
|\lambda_{a}|&<|\lambda_{c_{1}}|,\\
0<n_{c_{1}}-n_{a}-|J_{c_{1}}|k_{c_{1}}&<n_{a},\\
\alpha=n_{c_{1}}-n_{a}-|J_{c_{1}}|k_{c_{1}}&\leq\mathrm{max}\{\lambda_{a}\}.
\end{split}
\end{equation}
\par\bigskip

The representation $\Lambda'$ of the derived poset $\mathscr{P}'_{(a,b)}$ is given by the following formulas,
($2\leq j\leq i$),

\begin{equation}
\begin{split}
\Lambda'_{0} & =\Lambda_{0}, \\\notag
 n'_{a^{-}}&
=\alpha,\quad
\lambda'_{a^{-}}=((\alpha)^{1}),\\
n'_{c^{-}_{1}}&=n_{a},\hspace{1.5cm} \lambda'_{c^{-}_{1}}=\lambda_{a},\\
n'_{c_{i}^{-}} & =|\text{I}_{c_{i}}|k_{c_{i}},\hspace{1.0cm}
\lambda'_{c^{-}_{i}}= ((k_{c_{i}})^{|\mathrm{I}_{c_{i}}|}),\\
n'_{c_{i}^{+}} &=n_{c_{i}},\hspace{1.4cm} \lambda'_{c^{+}_{i}} =((\alpha)^{1}\lambda'_{c^{-}_{1}}(|\mathrm{I}_{c_{j}}|(k_{c_{j}}-k_{c_{j-1}}))^{\delta_{i}}k^{|J_{c_{i}}|}_{c_{i}}),\\
(n'_{x},\lambda'_{x})& =(h_{x},(h_{x}^{1})),\hspace{0.3cm}\text{for all}\hspace{0.1cm}x\in
b_{\vartriangle},\hspace{0.1cm}\text{if}\hspace{0.1cm}B^{\triangledown}\cap\{a^{-}\}\neq\varnothing,\\
h_{x}&=\mathrm{min}\{\alpha,\mathrm{min}\{n_{a}-n_{x}\mid x\in B\}\},\\
(n'_{x},\lambda'_{x})&=(n_{x},\lambda_{x}),\hspace{0.6cm}\text{for all}\hspace{0.1cm}x\in
b_{\vartriangle},\hspace{0.1cm}\text{if}\hspace{0.1cm}B^{\triangledown}\cap\{a^{-}\}=\varnothing.
\end{split}
\end{equation}
\begin{equation*}
\delta_{i}=\begin{cases}
1,&\text{if}\hspace{0.2cm} i\geq2, \\
0,&\text{otherwise}.
\end{cases}
\end{equation*}

(if $i=1$ then the formula for $\lambda'_{c^{+}_{1}}$ is given by $\lambda'_{c^{+}_{i}}$ by deleting the part
corresponding to $|\mathrm{I}_{c_{1}}|(k_{c_{j}}-k_{c_{j-1}})$, $n'_{c^{+}_{1}}=n_{c_{1}}$).\par\bigskip

Let $\mathscr{P}_{c_{i}}$ be the family of posets with a Hasse diagram of the form\par\smallskip

 \setlength{\unitlength}{1pt}
  \setlength{\unitlength}{1pt}
  \begin{centering}
\begin{picture}(95,95)

\put(24,30){\line(0,1){20}} \put(24,27){\circle{5}} \put(24,56){\line(0,1){20}} \put(24,53){\circle{5}}
\put(24,79){\circle{5}} \put(22,26){\line(-1,-1){20}} \put(0,4){\circle{5}} \put(0,-8){$a_{1}$}
\put(31,25){$c_{1}$} \put(31,50){$c_{2}$} \put(31,78){$c_{n}$} \put(54,2){$b$} \put(78,-22){$b_{1}$}
\put(24,-32){$\mathscr{P}_{c_{1}}$} \put(48,4){\circle{5}} \put(50,1){\line(1,-1){20}}
\put(26,26){\line(1,-1){20}}

\put(72,-21){\circle{5}}

\put(174,30){\line(0,1){20}} \put(174,27){\circle{5}} \put(174,56){\line(0,1){20}} \put(174,53){\circle{5}}
\put(174,79){\circle{5}} \put(172,51){\line(-1,-1){20}} \put(150,29){\circle{5}} \put(150,17){$a_{2}$}
\put(180,25){$c_{1}$} \put(180,50){$c_{2}$} \put(180,78){$c_{n}$} \put(204,2){$b$} \put(228,-22){$b_{1}$}
\put(174,-32){$\mathscr{P}_{c_{2}}$} \put(198,4){\circle{5}} \put(200,1){\line(1,-1){20}}
\put(176,26){\line(1,-1){20}}

\put(222,-21){\circle{5}} \put(174,-50){Fig. 3}

\put(202,25){$\dots$} \put(286,30){\line(0,1){20}} \put(286,27){\circle{5}} \put(286,56){\line(0,1){20}}
\put(286,53){\circle{5}} \put(286,79){\circle{5}} \put(283,77){\line(-1,-1){20}} \put(262,54){\circle{5}}
\put(262,42){$a_{n}$} \put(292,25){$c_{1}$} \put(292,50){$c_{2}$} \put(292,78){$c_{n}$} \put(316,2){$b$}
\put(340,-22){$b_{1}$} \put(286,-32){$\mathscr{P}_{c_{n}}$} \put(310,4){\circle{5}} \put(312,1){\line(1,-1){20}}
\put(288,26){\line(1,-1){20}}

\put(334,-21){\circle{5}}
\end{picture}
\end{centering}

\par\bigskip
\par\bigskip
\par\bigskip
\par\bigskip
\par\bigskip
\par\bigskip
We represent $\mathscr{P}_{c_{1}}$ over $\mathbb{N}$ in such a way that $\Lambda_{0}\supset\{t\in\mathbb{N}\mid
t=p^{n_{0}}_{r}\}$ (i.e., $\Lambda_{0}$ contains the set of the $n_{0}$-gonal numbers with $n_{0}\geq5$, fixed),

\begin{equation}\label{formulas 16}
\begin{split}
(n_{b_{1}},\lambda_{b_{1}})&=(n_{0}p^{n_{0}}_{1},((p^{n_{0}}_{1})^{n_{0}})),\\
(n_{b},\lambda_{b})&=(n_{0}p^{n_{0}}_{2},((p^{n_{0}}_{2})^{n_{0}})),\\
(n_{a_{1}},\lambda_{a_{1}})&=(p^{n_{0}}_{2}+p^{n_{0}}_{3},((p^{n_{0}}_{2})^{1}(p^{n_{0}}_{3})^{1})),\\
(n_{c_{i}},\lambda_{c_{i}})&=(n_{0}p^{n_{0}}_{i+2},((p^{n_{0}}_{i+2})^{n_{0}})),\hspace{0.2cm}i\geq1.
\end{split}
\end{equation}
\par\bigskip

For the corresponding representation of the poset $\mathscr{P}_{c_{i}}$, $2\leq i\leq n$, we choose the same
$\Lambda_{0}$, $n_{a_{i}}=p^{n_{0}}_{i+2}+p^{n_{0}}_{\iota}$,\quad
$\lambda_{a_{i}}=((p^{n_{0}}_{i+2})^{1}(p^{n_{0}}_{\iota})^{1})$, $2\leq \iota< i+2$, and leave the formulas
(\ref{formulas 16}) without changes for each $x\in b^{\triangledown}_{1}$.\par\bigskip

The next results were obtained in [14], with the help of the differentiation $\mathrm{L}$ as well as of the
representations of the posets $\mathscr{P}_{c_{i}}$ over $\mathbb{N}$ given above.

\begin{teor}
If $p^{3}_{0}=p^{3}_{-1}=0$, $p^{3}_{1}=1$, and $i,j\geq1$ then the entries of the matrices $R=(r_{(j,i}))$,
$S=(s_{(j,i)})$ and $T=(t_{(j,i)})$ satisfy the following identities
\begin{equation}
\begin{split}
24(r_{(j,i)}-2p^{5}_{i+1})+3&=(6i+5)^{2}+(6j+11)^{2}+(12j+29)^{2},\\
8(s_{(j,i)}-3p^{6}_{i+1})+3&=(4i+3)^{2}+(4j+7)^{2}+(8j+19)^{2},\\
40(t_{(j,i)}-4p^{7}_{i+1})+27&=(10i+7)^{2}+(10j+17)^{2}+(20j+47)^{2}
\end{split}
\end{equation}
where
\begin{equation}
\begin{split}
r_{(j,i)}&=97+57(j-1)+15p^{3}_{j-2}+9p^{3}_{i-2}+21(i-1),\\
s_{(j,i)}&=130+75(j-1)+20p^{3}_{j-2}+16p^{3}_{i-2}+36(i-1),\\
t_{(j,i)}&=165+93(j-1)+25p^{3}_{j-2}+25p^{3}_{i-2}+55(i-1).\hspace{0.5cm}\text{\qed}
\end{split}
\end{equation}

\end{teor}

\begin{teor}

If $i,j\geq1$ then $x=-(4j+8)$, $y=-(2j+4)$, and $z=2j+i+8$ is a solution of the diophantine equation\par\bigskip

\begin{centering}

$x^{3}+y^{3}+2z^{3}=m_{(j,i)}-n_{(j,i)}$,\par\bigskip

\end{centering}
where\par\smallskip $m_{(j,i)}=508+(j+1)891+690p^{3}_{j}+
(507+306j+72p^{3}_{j-1})(i+1)+(144+36j)p^{3}_{i}+18\rho^{3}_{i-1}+198\rho^{3}_{j-1}$, and\par\bigskip

$n_{(j,i)}=4068+4193(j-1)+1167(j-2)(j-1)+89(j-3)(j-2)(j-1)+(397+150(j-1)+12(j-2)(j-1))(i-1)+6(5+j)(i-2)(i-1)+(i-3)(i-2)(i-1)$.
That is,\par\bigskip

\begin{centering}

$M=\left[
  \begin{array}{ccccc}
    4786 & 5977 & 7384&\dots   \\
    8047 & 9688 &11581&\dots  \\
    12538 & 14701 & 17152&\dots  \\
     18457&21214&24295&\dots\\
     \vdots&\vdots&\vdots&\vdots
  \end{array}\right],\quad N=\left[
               \begin{array}{ccccc}
                 4068 & 4465 & 4934&\dots   \\
                 8261 & 8808 &9439&\dots  \\
                 14788 & 15509 & 16326&\dots  \\
24183&25102&26129&\dots\\
\vdots&\vdots&\vdots&\vdots
               \end{array}
             \right]$.\par\bigskip

  \end{centering}

\end{teor}

\textbf{Proof.} If $q_{k}$ is the $k$-th positive cube, each term $m_{(j,i)}$, $n_{(j,i)}$ can be obtained by
making the substitutions\par\bigskip
\begin{centering}
$\frac{\partial^{2}\lambda_{k}}{\partial q_{k+2}\partial
q_{k+2}}(q_{s},q_{2s+1})$,\quad$\frac{\partial^{2}\lambda_{ji}}{\partial q_{2j+i+8}\partial
q_{2j+i+8}}(q_{2j+4},q_{4j+8})$.\quad $s\geq 4$,\hspace{0.1cm} $2s-1\leq k$,\par\bigskip
\end{centering}

where $\lambda_{k},\lambda_{ji}$ are partitions such that
\begin{equation}
\begin{split}
\lambda_{k}&=((q_{k+2})^{5}),\\
\lambda_{ji}&=((q_{2j+7})^{1}(q_{j+3})^{1}(q_{2j+i+8})^{3}).\hspace{0.5cm}\text{\qed}
\end{split}
\end{equation}
\setcounter{corol}{4}

\begin{corol}

If\hspace{0.1cm} $i,j\geq1$, then $x=j+3$, $y=2j+7$\hspace{0.1cm} and \hspace{0.1cm}$z=2j+i+8$ is a solution of
the diophantine equation
\begin{equation}
\begin{split}
x^{3}+y^{3}+2z^{3}&=m_{(j,i)}-q_{2j+i+8}.\hspace{0.5cm}\text{\qed}
\end{split}
\end{equation}
\end{corol}

\begin{corol}\label{problema abierto}
If\hspace{0.1cm} $i,j\equiv 2\hspace{0.05cm}(\mathrm{mod}\hspace{0.05cm}3)$, $h\equiv
1\hspace{0.05cm}(\mathrm{mod}\hspace{0.05cm}3)$, $k=6+45l$, $m=2+9l$, $h,i,j,l\geq1$ then
\begin{equation}
\begin{split}
m_{(j,i)}-q_{2j+i+8}&\equiv-4\hspace{0.05cm}(\mathrm{mod}\hspace{0.05cm}9),\\
n_{(j,i)}-q_{4j+8}&\equiv-4\hspace{0.05cm}(\mathrm{mod}\hspace{0.05cm}9),\\
m_{(j,i)}-2q_{2j+i+8}+q_{2j+4}&\equiv -4\hspace{0.05cm}(\mathrm{mod}\hspace{0.05cm}9),\\
m_{(j,h)}-3q_{2j+h+8}+2q_{2j+h+k}&\equiv -4\hspace{0.05cm}(\mathrm{mod}\hspace{0.05cm}9),\\
m_{(j,i)}-3q_{2j+i+8}+2q_{2j+i+m}&\equiv -4\hspace{0.05cm}(\mathrm{mod}\hspace{0.05cm}9)\\
\end{split}
\end{equation}
and each number $m_{(j,i)}-q_{2j+i+8}$, $n_{(j,i)}-q_{4j+8}$, $m_{(j,i)}-2q_{2j+i+8}+q_{2j+4}$\hspace{0.1cm} and
\hspace{0.1cm} $m_{(j,i)}-3q_{2j+i+8}+2q_{2j+i+n}$, $8\neq n\geq1$ is the sum of four positive cubes.

\end{corol}

\par\bigskip
\textbf{3.1. The associated graph}\par\bigskip

Given a representation $\Lambda$, sometimes it is possible to associate it a suitable graph, having as set of
vertices the points of $\mathscr{P}$, and containing all information about parts and partitions of the numbers
$n_{x}$. In this case we must attach to each vertex of the graph, either a number $n_{x}$ given by the
representation or one part of a partition of some $n_{y}$ representing some $y\in\mathscr{P}$ such that $x\leq y$.
The theorem \ref{sigma} given more ahead shows the strength of this concept. To see this we must consider the
(infinite) poset $\mathscr{M}$ with Hasse diagram of the form

\par\bigskip
\par\bigskip
  \begin{centering}
\begin{picture}(100,100)

  \multiput(40,-80)(30,30){4}{\circle{5}}

\multiput(40,-20)(30,30){3}{\circle{5}}

\multiput(40,40)(30,30){2}{\circle{5}}

\multiput(40,100)(30,30){1}{\circle{5}}

\multiput(43,-77)(30,30){3}{\line(1,1){25}} \multiput(43,-17)(30,30){2}{\line(1,1){25}}
\multiput(43,43)(30,30){1}{\line(1,1){25}} \multiput(43,-23)(30,30){3}{\line(1,-1){25}}

\multiput(43,37)(30,30){2}{\line(1,-1){25}} \multiput(43,97)(30,30){1}{\line(1,-1){25}}
\multiput(46,-85)(30,30){1}{\tiny{$v_{11}$}} \multiput(76,-55)(30,30){1}{\tiny{$v_{21}$}}
\multiput(106,-25)(30,30){1}{\tiny{$v_{31}$}}
\multiput(136,5)(30,30){1}{\tiny{$v_{41}$}}\multiput(137,15)(5,5){3}{\tiny{$\cdot$}}
\multiput(30,-112)(30,30){1}{Fig. 4}

\multiput(0,5)(30,30){1}{$\mathscr{M}=$}

\multiput(76,5)(30,30){1}{\tiny{$v_{32}$}} \multiput(46,-25)(30,30){1}{\tiny{$v_{22}$}}
\multiput(106,35)(30,30){1}{\tiny{$v_{42}$}}\multiput(77,75)(5,5){3}{\tiny{$\cdot$}}
\multiput(46,35)(30,30){1}{\tiny{$v_{33}$}}\multiput(107,45)(5,5){3}{\tiny{$\cdot$}}
\multiput(47,105)(5,5){3}{\tiny{$\cdot$}}\multiput(27,106)(5,5){1}{\tiny{$\ddots$}}
\multiput(46,95)(30,30){1}{\tiny{$v_{44}$}} \multiput(76,65)(30,30){1}{\tiny{$v_{43}$}}

\end{picture}

\par\bigskip
\end{centering}
\par\bigskip
\par\bigskip
\par\bigskip
\par\bigskip
\par\bigskip
\par\bigskip
\par\bigskip
\par\bigskip
\par\bigskip
\par\bigskip

and a representation $\Delta$ of $\mathscr{M}$ over $\mathbb{N}$ such that, $\Delta_{0}$ contains all sum of the
form $\mathscr{O}_{i}+\mathscr{O}_{j}+\mathscr{O}_{k}$, (where $\mathscr{O}_{l}$ is the $l$-th octahedral number)
and for $i,j\geq1$, the pair
$(n_{ij},\lambda_{ij})=(\mathscr{O}_{ijk_{0}},(\mathscr{O}_{i})^{1}(\mathscr{O}_{j})^{1}(\mathscr{O}_{k_{0}})^{1})$
represents $v_{ij}\in\mathscr{P}$, where
$\mathscr{O}_{ijk_{0}}=\mathscr{O}_{i}+\mathscr{O}_{j}+\mathscr{O}_{k_{0}}$ and $k_{0}\geq1$ is a fixed number. It
is concluded the theorem \ref{sigma} due to that $\Delta$ induces the following partitions.\par\bigskip

If $n=\mathscr{O}_{ijk}$, then we note $P^{ijk}_{\mathscr{O}}(n)$ the number of partitions (of type $\mathscr{O}$)
of $n$ with the form
\begin{equation}\label{parttion O}
\begin{split}
n&=\mathscr{O}_{rsk}+z_{1}+\dots+z_{t},
\end{split}
\end{equation}

where if $z_{2j-1}\neq0$ for $1\leq j<\frac{t+1}{2}$ then $z_{2j}\neq0$, $z_{t}\neq0$, and\par\bigskip
\begin{centering}

 $z_{2j-1}=p^{4}_{k_{j}}$,\quad$z_{2j}=p^{4}_{(k_{j}+1)}$,\quad$|\{z_{m}\mid1\leq m\leq t\}|\leq2(j+i-2)$\par\bigskip

\end{centering}

for some $k_{j}$. Further, for each partial sum of the form\par\bigskip

$\mathscr{O}_{rsk}+z_{1}+\dots+z_{h}$, \quad$z_{h-1}=p^{4}_{k_{(h-1)}}$,\quad $z_{h}=p^{4}_{k_{(h-1)}+1}$,\quad
$h\leq t$\par\bigskip there exists $\mathscr{O}_{rsk}\leq\mathscr{O}_{pqk}\leq\mathscr{O}_{ijk}$ satisfying the
condition\par\bigskip
\begin{centering}

 $\mathscr{O}_{pqk}=\mathscr{O}_{rsk}+z_{1}+\dots+z_{h}$,\quad $r\leq p\leq i$,\quad $s\leq q\leq j$.\par\bigskip

\end{centering}

We denote now, $P^{iik+\sigma}_{5}(n)$ the number of partitions (of type $\sigma$) of $n=\mathscr{O}_{iik}$ with
the shape
\begin{equation}\label{partition P}
\begin{split}
n&=p^{5}_{rsk}+\sigma_{k}(i-1)+y_{1}+\dots+y_{t},
\end{split}
\end{equation}

where $p^{5}_{ijk}=p^{5}_{i}+p^{5}_{j}+p^{5}_{k}$, (i.e., is a sum of three pentagonal numbers)

\par\bigskip

\begin{centering}

$\sigma_{k}(i)=\underset{m=0}{\overset{k-1}{\sum}}p^{3}_{2m-1}+12i+18p^{3}_{i-3}-10$,\quad $p^{3}_{k}=0$ if
$k\leq0$. \par\bigskip

\end{centering}

If $y_{j}\neq0$ then $y_{j}=1+3\nu$, for some $\nu\geq1$,\quad $|\{y_{j}\mid1\leq j\leq t\}|\leq 2(i-1)$, and for
each partial sum of the form $p^{5}_{rsk}+\sigma_{k}(i-1)+y_{1}+\dots+y_{h}$, there exists $p^{5}_{rsk}\leq
p^{5}_{pqk}\leq p^{5}_{iik}$ satisfying the condition\par\bigskip

\begin{centering}

 $p^{5}_{pqk}+\sigma_{k}(i-1)=p^{5}_{rsk}+\sigma_{k}(i-1)+y_{1}+\dots+y_{h}$.\par\bigskip
\end{centering}

For example the partition $3+1^{2}+2^{2}+2^{2}+3^{2}+1^{2}+2^{2}+2^{2}+3^{2}$ is a partition of type $\mathscr{O}$
for $39=\mathscr{O}_{331}$, and $3+14+4+7+4+7$ is a partition of type $\sigma$.

\setcounter{teor}{5}

\begin{teor}\label{sigma}

$P^{iik}_{\mathscr{O}}(n)=P^{iik+\sigma}_{5}(n)$.

\end{teor}

\textbf{Proof.} If $\mathscr{M}$ is the poset described above represented either with three octahedral or three
pentagonal numbers, and $\Gamma_{\mathscr{M}}$ is an associated graph with the shape\par\bigskip
\begin{centering}
$\Gamma_{\mathscr{M}}=\begin{array}{ccccccccccc}
  \tiny{v_{11}}&\circ &\tiny{v_{22}}  & \circ & \tiny{v_{33}} & \circ &\tiny{v_{44}}& \circ&\dots&  \\
  &\downarrow &\nearrow  &\downarrow  &\nearrow  & &\nearrow  &  \\
 \tiny{v_{21}} &\circ &  & \circ &  & \circ &  &  \\
  &\downarrow &\nearrow &\downarrow &\nearrow  &\vdots  &  & \\
\tiny{v_{31}}& \circ&  & \circ & &  &  &  \\
  &\downarrow &\nearrow &  & &  & &  \\
\tiny{v_{41}}& \circ &  &\vdots  &  &  &  &&\\
&\vdots&
\end{array}$
\par\bigskip
\end{centering}

then the conclusion is obtained taking into account that $P^{iik}_{\mathscr{O}}(n)$ is the number of paths in
$\Gamma_{\mathscr{M}}$, with initial vertex of the form $v_{rs}$, $r\leq i$, $s\leq i$, and final vertex $v_{ii}$
(counting $v_{ii}$ as a path). This fact, because any partition of type $\mathscr{O}$ ($\sigma$) can be obtained
of a path $P=v_{rs}\|v_{ii}\in\Gamma_{\mathscr{M}}$, by attaching the number $\mathscr{O}_{rsk}$
($p^{5}_{rsk}+\sigma_{k}(i-1)$) to the vertex $v_{rs}$ and a suitable number of the form $z_{j}=p^{4}_{k_{j}}$
($y_{j}=1+3\nu$) to the rest of the vertices $v_{ij}\in P$. \hspace{0.5cm}\text{\qed}

\par\bigskip

The proof of the corollary \ref{partition pentagonal} given below, can be obtained fixing $q_{k}$ the $k$-th
positive cube and attaching to each $v_{ij}\in\mathscr{M}$ the pair
$(n_{ij},\lambda_{ij})=(\mathscr{Q}_{ijk},(q_{i})^{1}(q_{j})^{1}(q_{k})^{2})$, where
$\mathscr{Q}_{ijk}=q_{i}+q_{j}+2q_{k}$.

\setcounter{corol}{6}

\begin{corol}\label{partition pentagonal}

$P^{iik}_{\mathscr{Q}}(n)=P^{iik+\tau}_{\mathscr{O}}(n)$.\hspace{0.5cm}\text{\qed}

\end{corol}

Where $P^{iik}_{\mathscr{Q}}(n)$ is the number of partitions (of type $\mathscr{Q}$) of $n=\mathscr{Q}_{iik}$ with
the form,

\begin{equation}\label{partition Q}
\begin{split}
n&=\mathscr{Q}_{rsk}+y_{a_{1}}+\dots+y_{a_{t}},\quad a_{l}<a_{l+1},\quad l\geq1.
\end{split}
\end{equation}

If $y_{a_{j}}\neq0$ then $y_{a_{j}}=\frac{p^{5}_{2\nu+1}+\nu+1}{2}$, for some $\nu\geq1$,\quad
$|\{y_{a_{j}}\mid1\leq j\leq t\}|\leq 2(i-1)$, and for each partial sum of the form
$\mathscr{Q}_{rsk}+y_{a_{1}}+\dots+y_{a_{h}}$, there exists $\mathscr{Q}_{rsk}\leq \mathscr{Q}_{pqk}\leq
\mathscr{Q}_{iik}$ satisfying the condition\par\bigskip

\begin{centering}
 $\mathscr{Q}_{pqk}=\mathscr{Q}_{rsk}+y_{a_{1}}+\dots+y_{a_{h}}$, \quad$r\leq p\leq i$,\quad $s\leq q\leq i$.\par\bigskip

\end{centering}

$\tau_{k}(i)=\rho_{2k-1}+4(i-1)+8p^{3}_{i-2}+4p^{3}_{i-3}$, and $P^{iik+\tau}_{\mathscr{O}}(n)$ is the number of
partitions of type $\tau$ for $n=\mathscr{Q}_{iik}$ (i.e., partitions of type $\mathscr{O}$ for
$\mathscr{O}_{iik}$ with the additional term $\tau_{k}(i)$ such that the sum
$n=\mathscr{Q}_{iik}=\mathscr{O}_{rsk}+\tau_{k}(i)+z_{1}+\dots+z_{t}$ is defined as the sum
$p^{5}_{rsk}+\sigma_{k}(i-1)+y_{1}+\dots+y_{t}$ of three octahedral numbers).\par\bigskip

For example
$4+\frac{p^{5}_{3}+2}{2}+\frac{p^{5}_{5}+3}{2}+\frac{p^{5}_{3}+2}{2}+\frac{p^{5}_{5}+3}{2}=4+7+19+7+19$ is a
partition of type $\mathscr{Q}$ for $Q_{331}=56$ and
$3+17+(1^{2}+2^{2})+(2^{2}+3^{2})+(1^{2}+2^{2})+(2^{2}+3^{2})$ is a partition of type $\tau$.\par\bigskip

The following list shows the partitions of type $\mathscr{O}$ for $39$.\par\bigskip

$39$,\par\bigskip
$13+2^{2}+3^{2}+2^{2}+3^{2}$,\par\bigskip$26+2^{2}+3^{2}$,\par\bigskip$21+1^{2}+2^{2}+2^{2}+3^{2}$,\par\bigskip
$8+2^{2}+3^{2}+1^{2}+2^{2}+2^{2}+3^{2}$,\par\bigskip$8+1^{2}+2^{2}+2^{2}+3^{2}+2^{2}+3^{2}$,\par\bigskip
$3+1^{2}+2^{2}+1^{2}+2^{2}+2^{2}+3^{2}+2^{2}+3^{2}$,\par\bigskip
$3+1^{2}+2^{2}+2^{2}+3^{2}+1^{2}+2^{2}+2^{2}+3^{2}$.\par\bigskip

If
$\xi_{k}(i,j)=p^{6}_{k-1}+\underset{m=1}{\overset{j}{\sum}}p^{6}_{m-1}+\underset{h=j}{\overset{i}{\sum}}p^{6}_{h}$,\quad
$j\leq i$,\quad $k\geq1$ fixed, then a partition of type $\xi$ for $n=\mathscr{O}_{ijk}$, it is obtained of a
partition $\sigma$ for $\mathscr{O}_{iik}$, by replacing in (\ref{partition P}) the term $\xi_{k}(i,j)$ for
$\sigma_{k}(i-1)$ and leaving the others terms and conditions without changes, taking into account that in this
case $|\{y_{s}\mid1\leq s\leq t\}|\leq (j+i-2)$, and $r\leq p\leq i$,\quad $s\leq q\leq j$. We note
$P^{ijk+\xi}_{5}(n)$ the number of partitions of type $\xi$ for $n=\mathscr{O}_{ijk}$. \setcounter{teor}{7}

\begin{teor}\label{type sigma j}

$P^{ijk}_{\mathscr{O}}(n)=P^{ijk+\xi}_{5}(n)$.
\end{teor}

\textbf{Proof.} It is enough to observe the number of paths in $\Gamma_{\mathscr{M}}$ (see theorem (\ref{sigma})),
with initial vertex $v_{rs}$ and final vertex $v_{ij}$,\quad $r\leq i$,\quad $s\leq j$.\hspace{0.5cm}\text{\qed}

\par\bigskip
Given a number with the shape $\mathscr{O}_{(j+s)jk}$, $s\geq0$, now we can dedicate our efforts in finding a
formula for its partitions of type $\mathscr{O}$ (For a partition of type $\mathscr{O}$, we assume
$|\{z_{m}\mid1\leq m\leq t\}|\leq2(2j+s-2)$). To do this, we give the following two definitions.\par\bigskip

Let $\Lambda$ be a representation over $\mathbb{N}$ for a poset $\mathscr{P}$, then the number\par\bigskip

\begin{centering}
$W_{\Lambda}(\mathscr{P})=\underset{x\in\mathscr{P}}{\sum}n_{x}$,\par\bigskip

\end{centering}

is the \textit{weight} of the representation $\Lambda$. For example, if $\mathscr{N}_{i}$ is a poset with Hasse
diagram of the form
\par\bigskip

 \setlength{\unitlength}{1pt}
  \setlength{\unitlength}{1pt}
  \begin{centering}
\begin{picture}(100,100)
\multiput(180,-30)(0,30){5}{\circle{5}} \multiput(180,-27)(0,30){4}{\line(0,1){24}}

\multiput(183,-29)(0,30){1}{\line(1,1){24}} \put(165,-59){Fig. 5} \multiput(209,-2)(0,30){1}{\circle{5}}
\put(165,-30){$c_{1}$}\put(165,0){$c_{2}$}\put(115,30){$\mathscr{N}_{i}=$}
\put(165,30){$c_{3}$}\put(165,40){$\vdots$}\put(165,60){$c_{i}$}\put(157,90){$c_{i+1}$}\multiput(214,2)(0,30){1}{$d_{i}$}
\end{picture}
\end{centering}
\par\bigskip

\par\bigskip
\par\bigskip
\par\bigskip
\par\bigskip
\par\bigskip
and representation over $\mathbb{N}$, $\Lambda$ given by the formulas

\begin{equation}
\begin{split}
(n_{c_{j}},\lambda_{c_{j}})&=(j,(j)^{1}),\hspace{0.2cm}\text{if}\hspace{0.2cm}1\leq j\leq i,\\
(n_{c_{i+1}},\lambda_{c_{i}+1})&=(i,(i)^{1}),\\
(n_{d_{i}},\lambda_{d_{i}})&=(i,(i)^{1})
\end{split}
\end{equation}

then $W_{\Lambda}(\mathscr{N}_{i})=p^{3}_{i}+2i$.\par\bigskip

If $X_{ij}=P^{ijk}_{\mathscr{O}}(n)$ is the number of partitions of type $\mathscr{O}$ for $\mathscr{O}_{ijk}$,
where $k$ is a fixed number and $j\geq2$ then we write

\begin{equation}
\begin{split}
\delta_{(j+t)j}&=\underset{m=0}{\overset{t}{\sum}}(-1)^{t+m}\binom{t}{m}X_{(j+m)j},
\end{split}
\end{equation}

where $\delta_{jj}=X_{jj}$ if $j\geq2$.\par\bigskip

For example, $\delta_{22}=3$,\quad$\delta_{33}=8$,\quad$\delta_{44}=22$. \par\bigskip

It is easy to see that if $j=1$, $s\geq0$, then $X_{(s+1)1}=s+1$. While for $j\geq2$ we have the next result
:\setcounter{corol}{8}
\begin{corol}
If $s\geq1$ then $P^{(j+s-1)jk}_{\mathscr{O}}(n)=\underset{h=0}{\overset{j}{\sum}}b_{h}(s)\delta_{(h+j)j}$,
where\par\bigskip

$b_{0}(s)=1$, $b_{1}(s)=s-1$, and $b_{2}(s)=p^{3}_{s-2}$ if $j=2$.\par\bigskip $b_{0}(s)=1$, $b_{1}(s)=s-1$,
$b_{h}(s)=p^{3}_{s-h}$, $2\leq h< j$ and $b_{j}(s)=\rho_{s-j}$ if $j\geq3$.

\end{corol}

\textbf{Proof.} If $\Delta$ is the representation over $\mathbb{N}$ for the poset $\mathscr{M}$, and
$\Gamma_{\mathscr{M}}$ is the corresponding associated graph as before, then the relations
\begin{equation}
\begin{split}
P^{(j+s)jk}_{\mathscr{O}}(n)&=P^{(j+s-1)jk}_{\mathscr{O}}(n)+P^{(j+s)(j-1)k}_{\mathscr{O}}(n)+1,\hspace{0.2cm}s\geq0,\\
P^{(s+2)2k}_{\mathscr{O}}(n)&=W_{\Lambda}(\mathscr{N}_{s+1}),\\
P^{jjk}_{\mathscr{O}}(n)&=P^{(j)(j-1)k}_{\mathscr{O}}(n)+1
\end{split}
\end{equation}

observed in the representation of $\Gamma_{\mathscr{M}}$, allow us to obtain both the terms $X_{(j+s)j}$, $0\leq
s\leq j+2$, and the finite sequence $d_{0}=X_{jj}$,\hspace{0.2cm}$d_{s}=\{X_{(j+s+1)j}-X_{(j+s)j}\mid0\leq s\leq
j+1\}$ which give raise to the terms $\delta_{(j+s)j}$, $b_{h}(s)$ of
$P^{(j+s)jk}_{\mathscr{O}}(n)=X_{(j+s)j}$.\hspace{0.5cm}\text{\qed}

\setcounter{Nota}{9}

\begin{Nota}
$P^{(s+2)2k}_{\mathscr{O}}(n)=W_{\Lambda}(\mathscr{N}_{s+1})=3+4s+p^{3}_{s-1}$,\quad $s\geq0$.

\end{Nota}

The following list presents the values of $P^{ijk}_{\mathscr{O}}(n)$, for $2\leq i\leq11$, $2\leq j\leq 5$ and a
fixed index $k$.\par\bigskip

  $P^{22k}_{\mathscr{O}}(n)=3$,\quad $P^{33k}_{\mathscr{O}}(n)= 8$, \quad $P^{44k}_{\mathscr{O}}(n)=22$,\quad$P^{55k}_{\mathscr{O}}(n)=
  64$.\par\bigskip
  $P^{32k}_{\mathscr{O}}(n)=7$,\quad  $P^{43k}_{\mathscr{O}}(n)=21$,\quad$P^{54k}_{\mathscr{O}}(n)= 63$,\quad$P^{65k}_{\mathscr{O}}(n)=
  195$.\par\bigskip
  $P^{42k}_{\mathscr{O}}(n)=12$,\quad
  $P^{53k}_{\mathscr{O}}(n)=40$,\quad$P^{64k}_{\mathscr{O}}(n)=130$,\quad$P^{75k}_{\mathscr{O}}(n)=427$.\par\bigskip
  $P^{52k}_{\mathscr{O}}(n)=18,\quad P^{63k}_{\mathscr{O}}(n)=  66,\quad P^{74k}_{\mathscr{O}}(n)=  231,\quad P^{85k}_{\mathscr{O}}(n)=
  803$.\par\bigskip
  $P^{62k}_{\mathscr{O}}(n)=25,\quad P^{73k}_{\mathscr{O}}(n)=  100,\quad P^{84k}_{\mathscr{O}}(n)=  375,\quad P^{95k}_{\mathscr{O}}(n)=
  1376$.\par\bigskip
  $P^{72k}_{\mathscr{O}}(n)=33,\quad P^{83k}_{\mathscr{O}}(n)=  143,\quad P^{94k}_{\mathscr{O}}(n)=  572,\quad
  P^{(10)5k}_{\mathscr{O}}(n)=2210$.
  \par\bigskip
$P^{82k}_{\mathscr{O}}(n)=42,\quad P^{93k}_{\mathscr{O}}(n)=  196,\quad P^{(10)4k}_{\mathscr{O}}(n)=833,\quad
P^{(11)5k}_{\mathscr{O}}(n)=
  3381$.\par\bigskip
  $P^{92k}_{\mathscr{O}}(n)=52,\quad P^{(10)3k}_{\mathscr{O}}(n)=260,\quad P^{(11)4k}_{\mathscr{O}}(n)=  1170$.\par\bigskip

Now, we consider an infinite sum of infinite chains pairwise incomparable $\mathscr{R}$ in such a way that
$\mathscr{R}=\underset{i=0}{\overset{\infty}{\sum}}C_{i}$, where $C_{j}$ is a chain such that
$C_{j}=v_{0j}<v_{1j}<v_{2j}<\dots$.\par\bigskip It is defined a representation over $\mathbb{N}$ for
$\mathscr{R}$, by fixing a number $n\geq3$ and assigning to each $v_{ij}$ the pair
$(n_{ij},\lambda_{ij})=(3+(n-2)i+(n-1)j,(3+(n-2)i+(n-1)j)^{1})$, we note $\mathscr{R}_{n}$ this representation,
and write $v_{ij}\in\mathscr{R}_{n}$ whenever it is assigned the number $n_{ij}=3+(n-2)i+(n-1)j$ to the point
$v_{ij}\in\mathscr{R}$ in this representation. Fig. 6 below suggests the Hasse diagram for this poset with its
associated graph $\Gamma_{p^{n}_{k}}$ which attaches to each vertex $v_{ij}$ the number $3+(n-2)i+(n-1)j$.
\par\bigskip
\par\bigskip
 \setlength{\unitlength}{1pt}
  \setlength{\unitlength}{1pt}
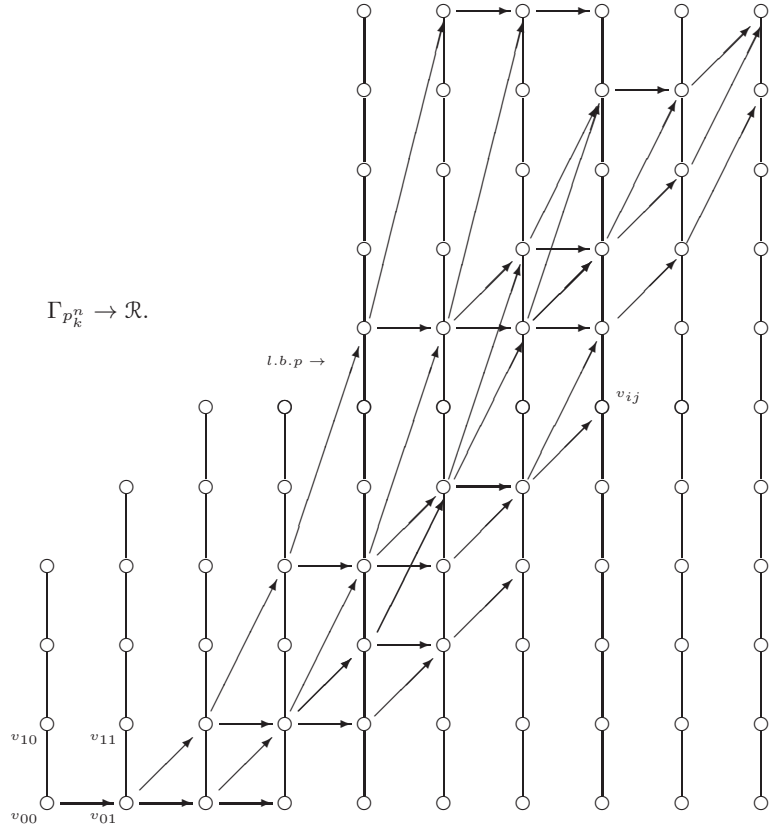
  \begin{centering}
\begin{picture}(315,315)

\multiput(40,110)(30,0){10}{\circle{5}} \multiput(40,80)(30,0){10}{\circle{5}}
\multiput(40,50)(30,0){10}{\circle{5}} \multiput(40,20)(30,0){10}{\circle{5}}
\multiput(70,140)(30,0){9}{\circle{5}} \multiput(100,170)(30,0){8}{\circle{5}}
\multiput(130,170)(30,0){6}{\circle{5}}

\multiput(160,200)(30,0){6}{\circle{5}} \multiput(160,230)(30,0){6}{\circle{5}}

\multiput(160,320)(30,0){6}{\circle{5}} \multiput(160,290)(30,0){6}{\circle{5}}
\multiput(160,260)(30,0){6}{\circle{5}}
 \put(40,204){$\Gamma_{p^{n}_{k}}\rightarrow\mathscr{R}$.}
\put(45,20){\vector(1,0){20}} \put(75,20){\vector(1,0){20}}\put(105,20){\vector(1,0){20}}
\put(75,25){\vector(1,1){20}} \put(102,55){\vector(1,2){25}} \put(135,55){\vector(1,1){20}}
\put(165,85){\vector(1,2){25}} \put(132,115){\vector(1,3){26}} \put(105,25){\vector(1,1){20}}
\put(132,55){\vector(1,2){25}} \put(162,115){\vector(1,3){26}} \put(105,50){\vector(1,0){20}}
\put(135,50){\vector(1,0){20}} \put(164,53){\vector(1,1){22}} \put(194,83){\vector(1,1){22}}
\put(194,113){\vector(1,1){22}} \put(194,203){\vector(1,1){22}}  \put(194,143){\vector(1,2){26}}
\put(192,143){\vector(1,3){27}} \put(165,80){\vector(1,0){20}} \put(195,140){\vector(1,0){20}}
\put(195,200){\vector(1,0){20}}

\put(224,143){\vector(1,1){22}} \put(222,143){\vector(1,2){26}} \put(222,204){\vector(1,3){27}}

 \put(225,230){\vector(1,0){20}} \put(225,200){\vector(1,0){20}} \put(224,203){\vector(1,1){22}}  \put(223,234){\vector(1,2){25}}

\put(162,204){\vector(1,4){28}} \put(192,204){\vector(1,4){28}}

\put(165,200){\vector(1,0){20}} \put(255,290){\vector(1,0){20}} \put(284,292){\vector(1,1){22}}
\put(284,262){\vector(1,2){26}}

 \put(282,232){\vector(1,2){26}}
 \put(256,204){\vector(1,1){22}} \put(256,234){\vector(1,1){22}}
\put(252,234){\vector(1,2){26}} \put(135,110){\vector(1,0){20}} \put(165,110){\vector(1,0){20}}
\put(165,114){\vector(1,1){22}}
 \put(195,320){\vector(1,0){20}} \put(225,320){\vector(1,0){20}}
 \multiput(40,23)(0,30){3}{\line(0,1){24}}

\multiput(70,23)(0,30){4}{\line(0,1){24}}

\multiput(100,23)(0,30){5}{\line(0,1){24}} \multiput(130,23)(0,30){5}{\line(0,1){24}}
\multiput(160,23)(0,30){10}{\line(0,1){24}} \multiput(190,23)(0,30){10}{\line(0,1){24}}
\multiput(220,23)(0,30){10}{\line(0,1){24}} \multiput(250,23)(0,30){10}{\line(0,1){24}}
\multiput(280,23)(0,30){10}{\line(0,1){24}} \multiput(310,23)(0,30){10}{\line(0,1){24}}

\put(135,55){\vector(1,1){20}}  \put(165,85){\vector(1,2){25}}

\put(224,203){\vector(1,1){22}}

\put(255,174){\tiny{$v_{ij}$}}
\put(26,44){\tiny{$v_{10}$}}\put(56,44){\tiny{$v_{11}$}}\put(56,14){\tiny{$v_{01}$}}
\put(26,14){\tiny{$v_{00}$}}\multiput(40,23)(0,30){3}{\line(0,1){24}}\put(123,186){\tiny{$l.b.p\rightarrow$}}

\put(40,-10){$v_{ij}\rightarrow3+(n-2)i+(n-1)j$,\quad $i,j\geq0$, $n\geq3$.} \put(40,-30){$v_{p^{3}_{(i-1)i}}\in
l.b.p$}

 \put(165,-55){Fig. 6}

\par\bigskip
\end{picture}
\end{centering}
\par\bigskip
\par\bigskip
\par\bigskip
\par\bigskip
\par\bigskip
\par\bigskip
The representations of $\mathscr{R}$ and $\Gamma_{p^{n}_{k}}$, allow us to observe that for $n_{0}$ fixed, each
natural number $n\geq n_{0}$ represents at least one point in $\mathscr{R}$, and that the numbers
$n_{p^{3}_{i-1}i}$ representing the vertices $v_{p^{3}_{(i-1)i}}\in\mathscr{R}_{t}$ of the \textit{left boundary
path}, $l.b.p$ have the form $n_{p^{3}_{i-1}i}=p^{t}_{1}+p^{t}_{1}+p^{t}_{i+1}$, $t$ is an fixed index.
Furthermore if $i_{0}\geq0$ is a fixed number, and $v_{p^{3}_{(i_{0}-1)i_{0}}}\in l.b.p$ then
$p^{t}_{i_{0}+1}+p^{t}_{k+1}+p^{t}_{1}$, represents the vertex\footnote{$v_{rs}\in P,\hspace{0.2cm} (P')$ if and
only if there exists $k_{0}$,\hspace{0.2cm}($l_{0}$),\hspace{0.2cm} $0\leq k_{0}\leq j$, $(0\leq l_{0}\leq j)$
such that
$v_{rs}=v_{(p^{3}_{i_{0}-1}+p^{3}_{k_{0}-1})(i_{0}+k_{0})},\hspace{0.2cm}(v_{(p^{3}_{i_{0}-1}+p^{3}_{j-1}+p^{3}_{l_{0}-1})(i_{0}+j+l_{0})})$.}
$v_{(p^{3}_{i_{0}-1}+p^{3}_{k-1})(i_{0}+k)}$ in the path
$P=v_{p^{3}_{(i_{0}-1)i_{0}}}||v_{(p^{3}_{i_{0}-1}+p^{3}_{j-1})(i_{0}+j)}$, $0\leq k\leq j$, and
$p^{t}_{i_{0}+1}+p^{t}_{j+1}+p^{t}_{l+1}$, represents the vertex
$v_{(p^{3}_{i_{0}-1}+p^{3}_{j-1}+p^{3}_{l-1})(i_{0}+j+l)}\in
P'=v_{(p^{3}_{i_{0}-1}+p^{3}_{j-1})(i_{0}+j)}||v_{(p^{3}_{i_{0}-1}+2p^{3}_{j-1})(i_{0}+2j)}$, $0\leq l\leq j$
(note that the vertices of the form $v_{(2+7i+3k)(2(i+k)+3)}$, $i,k\geq0$, do not lie on any non-trivial path of
$\Gamma_{p^{n}_{k}}$). These facts prove the following theorem ; \setcounter{teor}{10}

\begin{teor}\label{component}
A number $m\in\mathbb{N}$ is the sum of three $n$-gonal numbers of positive rank if and only if $m$ represents a
vertex $v_{ij}\in\mathscr{R}_{n}$ in a non-trivial component of $\Gamma_{p^{n}_{k}}$.\hspace{0.5cm}\text{\qed}

\end{teor}

\par\bigskip
 \setlength{\unitlength}{1pt}
  \setlength{\unitlength}{1pt}
  \begin{centering}
\begin{picture}(315,315)

\multiput(40,110)(30,0){10}{\circle{5}} \multiput(40,80)(30,0){10}{\circle{5}}
\multiput(40,50)(30,0){10}{\circle{5}} \multiput(40,20)(30,0){10}{\circle{5}}
\multiput(70,140)(30,0){9}{\circle{5}} \multiput(100,170)(30,0){8}{\circle{5}}
\multiput(130,170)(30,0){6}{\circle{5}}

\multiput(160,200)(30,0){6}{\circle{5}} \multiput(160,230)(30,0){6}{\circle{5}}

\multiput(160,320)(30,0){6}{\circle{5}} \multiput(160,290)(30,0){6}{\circle{5}}
\multiput(160,260)(30,0){6}{\circle{5}}
 \put(40,204){$\Gamma_{p^{n}_{k}}\rightarrow\mathscr{R}$.}
\put(45,20){\vector(1,0){20}} \put(75,20){\vector(1,0){20}}\put(105,20){\vector(1,0){20}}
\put(75,25){\vector(1,1){20}} \put(102,55){\vector(1,2){25}} \put(135,55){\vector(1,1){20}}
\put(165,85){\vector(1,2){25}} \put(132,115){\vector(1,3){26}} \put(105,25){\vector(1,1){20}}
\put(132,55){\vector(1,2){25}} \put(162,115){\vector(1,3){26}} \put(105,50){\vector(1,0){20}}
\put(135,50){\vector(1,0){20}} \put(164,53){\vector(1,1){22}} \put(194,83){\vector(1,1){22}}
\put(194,113){\vector(1,1){22}} \put(194,203){\vector(1,1){22}}  \put(194,143){\vector(1,2){26}}
\put(192,143){\vector(1,3){27}} \put(165,80){\vector(1,0){20}} \put(195,140){\vector(1,0){20}}
\put(195,200){\vector(1,0){20}}

\put(224,143){\vector(1,1){22}} \put(222,143){\vector(1,2){26}} \put(222,204){\vector(1,3){27}}

 \put(225,230){\vector(1,0){20}} \put(225,200){\vector(1,0){20}} \put(224,203){\vector(1,1){22}}  \put(223,234){\vector(1,2){25}}

\put(162,204){\vector(1,4){28}} \put(192,204){\vector(1,4){28}}

\put(165,200){\vector(1,0){20}} \put(255,290){\vector(1,0){20}} \put(284,292){\vector(1,1){22}}
\put(284,262){\vector(1,2){26}}

 \put(282,232){\vector(1,2){26}}
 \put(256,204){\vector(1,1){22}} \put(256,234){\vector(1,1){22}}
\put(252,234){\vector(1,2){26}} \put(135,110){\vector(1,0){20}} \put(165,110){\vector(1,0){20}}
\put(165,114){\vector(1,1){22}}
 \put(195,320){\vector(1,0){20}} \put(225,320){\vector(1,0){20}}
 \multiput(40,23)(0,30){3}{\line(0,1){24}}

\multiput(70,23)(0,30){4}{\line(0,1){24}}

\multiput(100,23)(0,30){5}{\line(0,1){24}} \multiput(130,23)(0,30){5}{\line(0,1){24}}
\multiput(160,23)(0,30){10}{\line(0,1){24}} \multiput(190,23)(0,30){10}{\line(0,1){24}}
\multiput(220,23)(0,30){10}{\line(0,1){24}} \multiput(250,23)(0,30){10}{\line(0,1){24}}
\multiput(280,23)(0,30){10}{\line(0,1){24}} \multiput(310,23)(0,30){10}{\line(0,1){24}}

\put(165,-20){Fig. 7}

\put(40,-10){$v_{ij}\rightarrow3+3i+4j$}

\put(135,55){\vector(1,1){20}}  \put(165,85){\vector(1,2){25}}

\put(224,203){\vector(1,1){22}}

\put(192,8){\tiny$23$} \put(222,8){\tiny$27$} \put(253,8){\tiny$31$} \put(282,8){\tiny$35$} \put(316,8){\tiny$39$}
\put(102,8){\tiny$11$}\put(132,8){\tiny$15$} \put(162,8){\tiny$19$}\put(72,8){\tiny$7$} \put(42,8){\tiny$3$}

\put(192,38){\tiny$26$} \put(222,38){\tiny$30$} \put(253,38){\tiny$34$} \put(282,38){\tiny$38$}
\put(316,38){\tiny$42$} \put(102,38){\tiny$14$}\put(132,38){\tiny$18$}
\put(162,38){\tiny$22$}\put(72,38){\tiny$10$} \put(42,38){\tiny$6$}

\put(192,68){\tiny$29$} \put(222,68){\tiny$33$} \put(253,68){\tiny$37$} \put(282,68){\tiny$41$}
\put(316,68){\tiny$45$} \put(102,68){\tiny$17$}\put(132,68){\tiny$21$}
\put(162,68){\tiny$25$}\put(72,68){\tiny$13$} \put(42,68){\tiny$9$}

\put(192,98){\tiny$32$} \put(222,98){\tiny$36$} \put(253,98){\tiny$40$} \put(282,98){\tiny$44$}
\put(316,98){\tiny$48$} \put(102,98){\tiny$20$}\put(132,98){\tiny$24$}
\put(162,98){\tiny$28$}\put(72,98){\tiny$16$} \put(42,98){\tiny$12$}

\put(162,128){\tiny$31$}\put(192,128){\tiny$35$} \put(222,128){\tiny$39$} \put(253,128){\tiny$43$}
\put(282,128){\tiny$47$}\put(312,128){\tiny$51$}

\put(162,158){\tiny$34$}\put(192,158){\tiny$38$} \put(222,158){\tiny$42$} \put(253,158){\tiny$46$}
\put(282,158){\tiny$50$}\put(312,158){\tiny$54$}

\put(162,188){\tiny$37$}\put(192,188){\tiny$41$} \put(222,188){\tiny$45$} \put(253,188){\tiny$49$}
\put(282,188){\tiny$53$}\put(312,188){\tiny$57$}

\put(192,218){\tiny$44$} \put(222,218){\tiny$48$} \put(253,218){\tiny$52$}
\put(282,218){\tiny$56$}\put(312,218){\tiny$60$}

\put(192,248){\tiny$47$} \put(222,248){\tiny$51$} \put(253,248){\tiny$55$}
\put(282,248){\tiny$59$}\put(312,248){\tiny$63$}

\put(192,278){\tiny$50$} \put(222,278){\tiny$54$} \put(253,278){\tiny$58$}
\put(282,278){\tiny$62$}\put(312,278){\tiny$66$}

\put(192,308){\tiny$53$} \put(222,308){\tiny$57$} \put(253,308){\tiny$61$}
\put(282,308){\tiny$65$}\put(312,308){\tiny$69$} \put(252,172){\vector(1,2){26}}

\par\bigskip
\end{picture}
\end{centering}

\par\bigskip
\par\bigskip
\par\bigskip
\par\bigskip
In the example given above (Fig. 7) each $v_{ij}\in\mathscr{R}_{5}$. That is, $n_{ij}=3+3i+4j$, $i,j\geq0$, and
the number associated to each vertex in a non-trivial component of $\Gamma_{p^{n}_{k}}$, is the sum of three
pentagonal numbers of positive rank.

\setcounter{corol}{11}
\begin{corol}\label{three triangular}
The number $m$ representing  $v_{ij}\in\mathscr{R}_{3}$ is the sum of three triangular numbers $\geq1$ if and only
if the number $n$ representing $v_{ij}\in\mathscr{R}_{t}$ is the sum of three $t$-gonal numbers of positive rank.
\end{corol}

Note that, the structure of the $l.b.p$ gives the general form of the graph $\Gamma_{p^{n}_{k}}$. Thus, it is
enough to change the form of such left boundary path, to build different graphs of this type. For instance, we
note $\Gamma_{\mathscr{O}}$ a graph associated to $\mathscr{R}$ which $v_{00}$, $v_{01}$,
$v_{(p^{3}_{(2i-1)}+1)(i+1)}$, $i\geq1$ are the locations of the vertices of its left boundary path. In this case
the representation over $\mathbb{N}$ for $\mathscr{R}$ is such that the number $n_{ij}=3+4i+5j$, $i,j\geq0$
represents the point $v_{ij}$ (we note $\mathscr{R}_{\mathscr{O}}$ this representation, and write
$v_{ij}\in\mathscr{R}_{\mathscr{O}}$ in this situation). Hence if $v_{ij}$ is a vertex in a non trivial component
of $\Gamma_{\mathscr{O}}$ then $n_{ij}$ can be expressed as the sum of three octahedral numbers (in particular
$n_{(p^{3}_{(2i-1)}+1)(i+1)}=\mathscr{O}_{i+2}+2\mathscr{O}_{1}$, $i\geq1$).\par\bigskip

 Now, we represent the poset $\mathscr{R}$ in such a way that to each point $v_{ij}$, it is associated the number $n_{ij}=4+6i+7j$
($v_{ij}\in\mathscr{R}_{\mathscr{Q}}$, see $\mathscr{R}_{\mathscr{O}}$), thus the numbers
$n_{(p^{3}_{(2i-1)}+1)(i+1)}$\hspace{0.1cm}or\hspace{0.1cm} $(n_{(\rho_{s}-s)s})$, $i\geq1$, $s\geq2$,
representing the vertices in the $l.b.p$ can be presented in the form $q_{1}+q_{i+2}+2q_{1}$ (if $s=0$ then
$n_{(\rho_{s}-s)(s)}=n_{00}=4$, represents the vertex $v_{00}$, and $n_{01}=11$ represents the vertex $v_{01}$ if
$s=1$). These facts and the theorem \ref{component} prove the next theorems. \setcounter{teor}{12}
\begin{teor}\label{octahedral}
A number $m\in\mathbb{N}$ is the sum of three octahedral numbers if and only if $m$ represents a vertex
$v_{ij}\in\mathscr{R}_{\mathscr{O}}$ in a non-trivial component of
$\Gamma_{\mathscr{O}}$.\hspace{0.5cm}\text{\qed}

\end{teor}

\begin{teor}\label{cubes}
A number $m\in\mathbb{N}$ is the sum of four positive cubes with two of them equal if and only if $m$ represents a
vertex $v_{ij}\in\mathscr{R}_{\mathscr{Q}}$ in a non-trivial component of the graph
$\Gamma_{\mathscr{O}}$.\hspace{0.5cm}\text{\qed}

\end{teor}

\begin{teor}\label{cubes and octahedral numbers}
The number $m$ representing  $v_{ij}\in\mathscr{R}_{\mathscr{O}}$ is the sum of three octahedral numbers if and
only if the number $n$ representing $v_{ij}\in\mathscr{R}_{\mathscr{Q}}$ is the sum of four positive cubes with
two of them equal.\hspace{0.5cm}\text{\qed}

\end{teor}

The corollaries \ref{cuadrados}, \ref{sum}, and \ref{diophantine} of the theorems \ref{component},
\ref{octahedral}, and \ref{cubes} respectively are also consequences, of the Gauss's theorem for three triangular
numbers, and of the structures of the graphs $\Gamma_{p^{n}_{k}}$ and $\Gamma_{\mathscr{O}}$. Furthermore these
corollaries can be interpreted as algorithms which solve diophantine equations of the form,
$n=x^{2}+y^{2}+z^{2}$,\quad$n=x^{3}+y^{3}+2z^{3}$,\quad and\quad
$n=x(2x^{2}+1)/3+y(2y^{2}+1)/3+z(2z^{2}+1)/3$,\quad $x,y,z,n\geq0$.

\setcounter{corol}{15}

\begin{corol}\label{cuadrados}
If $n\in\mathbb{N}$ is the sum of three square numbers of positive rank then there exists
$v_{ij}\in\mathscr{R}_{4}$, and $k_{0}\geq0$ such that
\begin{equation}
\begin{split}
n_{ij}&=n,\\
i-3k_{0}&=p^{3}_{\alpha(k_{0})}+p^{3}_{\beta(k_{0})}+p^{3}_{\gamma(k_{0})},\notag\\
j+2k_{0}-3&=\alpha(k_{0})+\beta(k_{0})+\gamma(k_{0}),
\end{split}
\end{equation}

$\alpha(k_{0}),\beta(k_{0}),\gamma(k_{0})\geq-1$. Therefore\par\bigskip
\begin{centering}
$n=p^{4}_{\alpha(k_{0})+2}+p^{4}_{\beta(k_{0})+2}+p^{4}_{\gamma(k_{0})+2}$, \par\bigskip

\end{centering}
where
\begin{equation}
\begin{split}
n&=p^{4}_{\alpha(k_{0})+2}+2p^{4}_{1},\quad \text{if} \quad p^{3}_{\beta(k_{0})}=p^{3}_{\gamma(k_{0})}=0,\notag\\
n&=p^{4}_{\alpha(k_{0})+2}+p^{4}_{\beta(k_{0})+2}+1,\quad \text{if} \quad p^{3}_{\gamma(k_{0})}=0.
\end{split}
\end{equation}

If the poset $\mathscr{R}$ is represented in such a way that $n_{ij}=2i+j$, for each $i,j\geq0$ and a number
$n\in\mathbb{N}$ is not of the form $4^{k}(8m+7)$, $k,m\geq0$ then there are $i,j\geq0$, and $k_{0}\geq0$ such
that
\begin{equation}
\begin{split}
n_{ij}&=n,\\
i-k_{0}&=p^{3}_{\alpha(k_{0})}+p^{3}_{\beta(k_{0})}+p^{3}_{\gamma(k_{0})},\notag\\
j+2k_{0}-3&=\alpha(k_{0})+\beta(k_{0})+\gamma(k_{0}),
\end{split}
\end{equation}

$\alpha(k_{0}),\beta(k_{0}),\gamma(k_{0})\geq-1$. Therefore\par\bigskip
\begin{centering}
$n=p^{4}_{\alpha(k_{0})+1}+p^{4}_{\beta(k_{0})+1}+p^{4}_{\gamma(k_{0})+1}$, \par\bigskip

\end{centering}
where
\begin{equation}
\begin{split}
n&=p^{4}_{\alpha(k_{0})+1},\quad \text{if} \quad p^{3}_{\beta(k_{0})}=p^{3}_{\gamma(k_{0})}=0,\notag\\
n&=p^{4}_{\alpha(k_{0})+1}+p^{4}_{\beta(k_{0})+1},\quad \text{if} \quad
p^{3}_{\gamma(k_{0})}=0.\hspace{0.5cm}\text{\qed}
\end{split}
\end{equation}

\end{corol}

For example $n_{(24)8}=75\in\mathscr{R}_{4}$ in this case $i=24$, $j=8$ and $24-3(2)=6+6+6=3p^{3}_{3}$. Hence
$k_{0}=2$, and $9=8+2(2)-3$. Thus, $75=p^{4}_{3+2}+p^{4}_{3+2}+p^{4}_{3+2}=25+25+25$.

\begin{corol}\label{sum} If $n\in\mathbb{N}$ is the sum of three octahedral numbers then there exists $v_{ij}\in\mathscr{R}_{\mathscr{O}}$, and
$k_{0}\geq0$ such that
\begin{equation}
\begin{split}
n_{ij}&=n,\\
i-5k_{0}&=\rho_{\alpha(k_{0})}+\rho_{\beta(k_{0})}+\rho_{\gamma(k_{0})}-(\alpha(k_{0})+\beta(k_{0})+\gamma(k_{0})),\notag\\
j+4k_{0}&=\alpha(k_{0})+\beta(k_{0})+\gamma(k_{0})
\end{split}
\end{equation}

$\alpha(k_{0}),\beta(k_{0}),\gamma(k_{0})\geq-1$, and $\rho_{0}=\rho_{-1}=0$. Thus\par\bigskip
\begin{centering}
$n=\mathscr{O}_{\alpha(k_{0})+1}+\mathscr{O}_{\beta(k_{0})+1}+\mathscr{O}_{\gamma(k_{0})+1}$,\par\bigskip

\end{centering}

where
\begin{equation}
\begin{split}
n&=\mathscr{O}_{\alpha (k_{0})+1}+2,\quad \text{if} \quad \rho_{\beta(k_{0})}=\rho_{\gamma(k_{0})}=0,\notag\\
n&=\mathscr{O}_{\alpha(k_{0})+1}+\mathscr{O}_{\beta(k_{0})+1}+1,\quad \text{if} \quad
\rho_{\gamma(k_{0})}=0.\hspace{0.5cm}\text{\qed}
\end{split}
\end{equation}
\end{corol}

For example $n_{(7)(4)}=51\in\mathscr{R}_{\mathscr{O}}$, $7=(\rho_{3}-3)+(\rho_{1}-1)$, and $4=3+1$. Therefore,
$51=\mathscr{O}_{4}+\mathscr{O}_{2}+1=44+6+1$.

\begin{corol}\label{diophantine} If $n\in\mathbb{N}$ is the sum of four positive cubes with two of them equal then there exists $v_{ij}\in\mathscr{R}_{\mathscr{Q}}$, and
$k_{0}\geq0$ such that
\begin{equation}
\begin{split}
n_{ij}&=n,\\
i-7k_{0}&=\rho_{\alpha(k_{0})}+\rho_{\beta(k_{0})}+\rho_{\gamma(k_{0})}-(\alpha(k_{0})+\beta(k_{0})+\gamma(k_{0})),\notag\\
j+6k_{0}&=\alpha(k_{0})+\beta(k_{0})+\gamma(k_{0})
\end{split}
\end{equation}

thus\par\bigskip
\begin{centering}
$n=q_{\alpha(k_{0})+1}+q_{\beta(k_{0})+1}+2q_{\gamma(k_{0})+1}$,\par\bigskip

\end{centering}

where
\begin{equation}
\begin{split}
n&=q_{\alpha(k_{0})+1}+3,\quad \text{if} \quad \rho_{\beta(k_{0})}=\rho_{\gamma(k_{0})}=0,\notag\\
n&=q_{\alpha(k_{0})+1}+q_{\beta(k_{0})+1}+2,\quad \text{if} \quad \rho_{\gamma(k_{0})}=0.\hspace{0.5cm}\text{\qed}
\end{split}
\end{equation}
\end{corol}
Let us to illustrate the corollary \ref{diophantine}, by considering the vertex
$v_{(14)6}\in\mathscr{R}_{\mathscr{Q}}$. In this case $n_{(14)6}=130$, $k_{0}=0$, $\alpha(k_{0})=\beta(k_{0})=3$,
$\gamma(k_{0})=0$, and $14=2\rho_{3}-6$,\quad thus\quad$130=q_{4}+q_{4}+2q_{1}$.

\setcounter{Nota}{18}

\begin{Nota}
Since $k_{0}=0$ is one of such values of $k_{0}$ in the corollaries, \ref{cuadrados}, \ref{sum},
\ref{diophantine}, an interesting problem consists in finding all the values of $k_{0}$ satisfying the
requirements.

\end{Nota}

\textbf{References}

\par\bigskip

[1] G. Andrews, \textit{The Theory of Partitions}, Cambridge Univ. Press, Cambridge, 1991.

\par\smallskip

[2] T. Apostol, \textit{Introduction to Analytic Number Theory}, Springer, 1998, 1-337.

\par\smallskip

[3] B. Davey and H. Priestley, \textit{Introduction to Lattices and Order}, Cambridge Univer-\par\smallskip\hspace*{0.5cm}sity
Press, 2002.

\par\smallskip

[4] J.M. Deshouillers, F. Hennecart, B. Landreau, and I.G. Putu,  \textit{7373170279850},
\hspace*{0.5cm}Mathematics of computation \textbf{69}(2000), no. 229, 1735-1742.
\par\smallskip

[5] Duke. W, \textit{Some old and new results about quadratic forms}, Notices. Amer. Math. \hspace*{0.5cm}Soc
\textbf{44} (1997), 190-196.

\par\smallskip

[6] J. A. Ewell, \textit{On sums of triangular numbers and sums of squares}, A.M.M \textbf{99}
(1992),\par\smallskip\hspace*{0.3cm} no. 8, 752-757.

\par\smallskip
[7]  J. A. Ewell, \textit{A Trio of triangular number theorems}, A.M.M \textbf{105} (1998), no. 9, 848-\par\smallskip\hspace*{0.5cm}849.

\par\smallskip

[8] H. Farkas, \textit{Sums of squares and triangular numbers}, Journal of analytic
combina-\par\smallskip\hspace*{0.5cm}torics \textbf{1} (2006), no. 1, 1-11.

\par\smallskip

[9] R. Guy, \textit{Every number is expressible as a sum of how many polygonal numbers}, \hspace*{0.5cm}A.M.M
\textbf{101} (1994), 169-172.

\par\smallskip
[10] R. Guy, \textit{Nothing's new in number theory?}, A.M.M \textbf{105} (1998), no. 10, 951-954.
\par\smallskip
[11] R. Guy, \textit{Unsolved Problems in Number Theory}, 3rd ed., Springer-Verlag, New \hspace*{0.6cm}York, 2004.
\par\smallskip
[12] D.R. Heath-Brown, W.M. Lioen, and J.J. Te Riele, \textit{On solving the
diophantine}\par\smallskip\hspace*{0.6cm}\textit{equation $x^{3}+y^{3}+z^{3}=k$ on a vector computer}, Mathematics
of computation \textbf{61}\par\smallskip\hspace*{0.5cm} (1993), no. 203, 235-244.

\par\smallskip

[13] K. Koyama, \textit{On searching for solutions of the diophantine equation
$x^{3}+y^{3}+2z^{3}=$}\par\smallskip\hspace{0.6cm}\textit{$n$}, Mathematics of computation \textbf{69} (2000),
1735-1742.

 \par\smallskip
[14] A. Moreno, \textit{Descripci\'on categ\'orica de algunos algoritmos de diferenciaci\'on},
Diss.\par\smallskip\hspace*{0.5cm} Univ. Nacional. Colombia (2007), 1-157.

\par\smallskip

[15] B. Schr\"oder, \textit{Ordered Sets An Introduction}, Birkh\"auser, 2003, 1-391.

\par\smallskip
[16] D. Simson, \textit{Linear Representations of Partially Ordered Sets and Vector
Space}\par\smallskip\hspace*{0.6cm}\textit{Categories}, Gordon and Breach, London, 1992.

\end{document}